\documentclass[a4paper,10pt]{article} \usepackage{amsmath,amssymb,amsthm}
\usepackage{dsfont}
\usepackage[english]{babel}

\newtheorem{theo}{Theorem}[section] 
\newtheorem{defi}[theo]{Definition}
\newtheorem{lemm}[theo]{Lemma} 
\newtheorem{prop}[theo]{Proposition}
\newtheorem{coro}[theo]{Corollary}

\newcommand{\Na}{\mathbb N}                   

\newcommand{\Ra}{\mathbb R}                   

\newcommand{\Ca}{\mathbb C}                   

\newcommand{\scal}[1]{\langle #1 \rangle}

\newcommand{\finpreuve}{\hfill $\Box$}

\newcommand{\name}{$\underline{\qquad \qquad}$}

\begin{document}

\author{  Jean-Marc Bouclet \\
Institut de Math\'ematiques de Toulouse \\ 
118 route de Narbonne 
\\ 
F-31062 Toulouse Cedex 9 \\
e-mail: jean-marc.bouclet@math.univ-toulouse.fr}
\title{Absence of eigenvalue at the bottom of the continuous spectrum on asymptotically hyperbolic manifolds}

\maketitle

\begin{abstract} For a class of asymptotically hyperbolic manifolds, we show that the bottom of the continuous spectrum of the Laplace-Beltrami operator is not an eigenvalue.  Our approach  only uses properties of the operator near infinity and, in particular, does not require any global assumptions on the topology or the curvature, unlike previous papers on the same topic. 
\end{abstract}



\section{Introduction and main result}
The main purpose of this paper is to prove that on an asymptotically hyperbolic manifold $ (M,G) $ of dimension $n$ and for  perturbations $V$  of the Laplace-Beltrami operator which decay at infinity, we have the property
\begin{eqnarray}
 \big( - \Delta_G  + V \big) \psi = \frac{(n-1)^2}{4} \psi \qquad \mbox{and} \qquad \psi \in L^2 \qquad \Longrightarrow \qquad \psi = 0 . \label{resultformel} 
\end{eqnarray}
Here $V$ may be  a potential but also a second order differential operator, possibly with complex coefficients. In the case when $ - \Delta_G + V $ is selfadjoint, this means that the bottom of the essential spectrum is not an eigenvalue. That the essential spectrum is absolutely continuous follows for instance from \cite{FroeseHislop}. Let us recall that, in scattering theory, ruling out the presence of such an eigenvalue  is the first step in the study of the resolvent of the Laplacian at the bottom of the continuous spectrum, this being in turn  important  to analyze long time  properties  of dispersive  equations, such as the local energy decay for the wave equation. This is the main motivation of this paper.

 The property (\ref{resultformel}) has been considered in several papers, with $ V \equiv 0 $, under various curvature conditions  but  only  for (perturbations of) spherically symmetric metrics  \cite{Pinsky,Donnelly3} and simply connected manifolds \cite{Donnelly1,Donnelly2}, each time with the additional assumption that the curvature is globally negative. These are  restrictive conditions, in particular for the purpose of scattering theory  of asymptotically hyperbolic manifolds where one wishes to treat general scatterers, ie not to rely too much on what happens in a compact region. In a different direction from the previous results, Vasy and Wunsch \cite{VasyWunsch} have obtained an absence of eigenvalue condition on a very general class of negatively curved manifolds  for which no asymptotic behavior is needed (only a pinching condition on the curvature is required), but under the stronger assumption that $   \psi $ decays super exponentially. 
 Recently, Kumura \cite{Kumura2} has obtained fairly sharp conditions on the decay rate of the radial curvature to $-1$ leading to the absence or the existence of embedded eigenvalues in the bulk of the continuous spectrum, however his results do not  apply  to the bottom of the continuous spectrum $ (n-1)^2 / 4 $.

 Our main result in this paper is  that if the left hand side of (\ref{resultformel}) holds, then $ \psi $ decays super exponentially. Using the unique continuation result of Mazzeo \cite{Mazzeo}, this implies automatically that $ \psi $ must  vanish identically (if $ V$ is a second order operator then, for this unique continuation purpose,  its principal symbol has to be real, but lower order terms can have complex coefficients).  The  interest of our approach is that
  it depends only on data at infinity: the metric and the topology can be arbitrary in a compact set, no global condition on the curvature is needed  and we don't use any spherical symmetry nor simple connectedness. Furthermore, our proof of the super exponential decay of $ \psi $  is robust enough to handle non self-adjoint operators (here $V$ may have complex coefficients even in the second order terms) and does not use crucially the particular structure of the angular Laplacian. In addition, our decay condition  on $G$ to an exact warped product $ dr^2 + e^{2r} \overline{g} $ (see (\ref{decroissance})) is relatively weak (it is a short range pertubation in the Schr\"odinger operators terminology), and in any case much weaker than what happens in the conformally compact case where one has exponential decay.





Here are our assumptions on $ (M,G) $ and on the perturbations $V$.

We consider an  asymptotically hyperbolic Riemannian manifold $ (M,G) $ of the following form. We assume that $ M$ is smooth and that, for some smooth compact subset $ K \Subset M $ with boundary $ \partial K = S $ (with $S$ of dimension $ n - 1$), we have
\begin{eqnarray}
 (M \setminus K , G ) \ \ \mbox{is isometric to} \ \  \big( (R,+\infty) \times S , dr^2 + e^{2r} g (r) \big) , \label{metriquehyperbolique} 
\end{eqnarray} 
where   $ (g (r))_{r > R} $ is a family of Riemannian metrics on $ S $ depending smoothly on $r$ and which converges to a fixed metric $ \overline{g} $ as $r \rightarrow \infty $, in the sense that
\begin{eqnarray}
 || \partial_r^j \big( g (r) - \overline{g} \big) ||_{C^{\infty}(T^* S \otimes T^* S)} \leq C \scal{r}^{- \tau_0 - j} , \label{decroissance}
\end{eqnarray}
for each semi-norm $ || \cdot ||_{C^{\infty}(T^* S \otimes T^* S)} $ of the space of smooth sections of $ T^* S \otimes T^* S  $. Here $ \tau_0 $ is a positive real number that
$$ \tau_0 > 1 . $$
 We  note  that  metrics of the more general form
$$ G = a (x,\theta) dx^2 + e^x b_{j}(x,\theta) dx d \theta_j + e^{2x} h_{ij}(x,\theta, d \theta_i , d \theta_j) $$
where $ \big( \theta_1 , \ldots , \theta_{n-1} \big) $ are coordinates on $ S $ and $x $ a coordinate such that $ x \rightarrow \infty $ at infinity on $ M \setminus K $,    can be put under the normal
form $ d r^2 + e^{2r}g (r) $ under natural decay rates of the coefficients $ a$ to $1$ and $b_j$ to $0$ (see \cite{Bouclet} for more details).  

The perturbations $V$ which are allowed in (\ref{resultformel}) are as follows. First we assume that they are second order differential operators on $ M $ with smooth coefficients such that
\begin{eqnarray}
 - \Delta_G + V \ \ \mbox{is  elliptic on} \ \ M .  \label{adm1}
\end{eqnarray}
This implies first that if the left hand side of (\ref{resultformel}) holds then $ \psi  $ is smooth and more importantly that the operator satisfies the local unique continuation principle provided that its principal symbol is real (see for instance \cite[Theorem 17.2.6]{H}). 
Near infinity, ie on $ M \setminus K $, we assume that
\begin{eqnarray}
 V = e^{-2r} W (r) + a (r,\omega) \partial_r + b (r,\omega) , \qquad \omega \in S ,  \label{adm2}
\end{eqnarray}
where $ W (r) $ is, for each each $r$, a second order differential operator on $ S $ which reads in local coordinates
\begin{eqnarray}
 W (r) = \sum_{|\alpha| \leq 2} a_{\alpha}(r,\theta) D_{\theta}^{\alpha} ,
 \label{adm3} \end{eqnarray}
with
\begin{eqnarray}
 |\partial_r^j \partial_{\theta}^{\beta} a_{\alpha}(r,\theta) | \leq C_{j \beta} \scal{r}^{-\tau_0 -j} , \qquad  j \leq 2, \ \beta \in \Na^{n-1}, 
 \label{adm4}
\end{eqnarray}
locally uniformly with respect to $ \theta $. We also assume that
\begin{eqnarray}
|| \partial_r^j a (r,.) ||_{C^{\infty}(S)} \leq C \scal{r}^{-\tau_0-1-j} \ \ \ j = 0 , 1, \qquad || b(r,.)||_{C^{\infty}(S)} \leq C \scal{r}^{- \tau_0 -1} , \label{adm5}
\end{eqnarray}
for all semi-norms $ ||\cdot ||_{C^{\infty}(S)} $ of $ C^{\infty}(S) $. For simplicity we also assume that $ a_{\alpha}$, $a$ and $b$ are smooth  but no bound on higher order derivatives in $r$ will be used.
Note also that these coefficients can be complex valued; furthermore, we do not require $ V$ to be symmetric with respect to $ d {\rm vol}_G $, the Riemannian volume density. An operator $V$ satisfying (\ref{adm1}), (\ref{adm2}), (\ref{adm3}), (\ref{adm4}) and (\ref{adm5}) will be called an {\it admissible perturbation}.

\begin{theo} \label{vraitheoreme} Let $ (M,G) $ be a connected asymptotically hyperbolic manifold of dimension $n$ and  $ V $ be an admissible perturbation. If $ \psi \in L^2 (M, d {\rm vol}_G) $ satisfies
\begin{eqnarray}
 \big(- \Delta_G + V \big) \psi = \frac{(n-1)^2}{4} \psi , \label{equationvaleurpropre}
\end{eqnarray} 
then on $ M \setminus K $ we have, for all $ C > 0 $,
\begin{eqnarray}
 e^{Cr} \psi \in L^2 , \qquad \partial_r \big( e^{Cr} \psi \big) \in L^2 , \qquad  (-\Delta_G + V) \big( e^{Cr} \psi \big) \in L^2  . \label{superexponentielfinal}
\end{eqnarray}
\end{theo}


In (\ref{superexponentielfinal}), $ L^2 $ stands for $ L^2 (M \setminus K, d {\rm vol}_G) $, since $r$ is only defined on $ M \setminus K $, but this is sufficient for we are only interested in the behavior near infinity and we know that $ \psi $ is smooth. The super exponential decay (\ref{superexponentielfinal}) and the result of Mazzeo \cite[Corollary (11)]{Mazzeo} on Carleman estimates and unique continuation lead to

\begin{coro} \label{corollairevp} Let $ V $ be an admissible perturbation, with real principal symbol if $V$ is a second order  operator\footnote{If $V$ is only a first order operator, then no additional condition is required}. Then (\ref{resultformel}) holds true.
In particular, for $V \equiv 0$, $ (n-1)^2 / 4 $ is not an eigenvalue of $ - \Delta_G $.
\end{coro}

Note that even if we show that $ (n-1)^2/4 $  is not an eigenvalue of $ - \Delta_G $, we do not exclude that there can be eigenvalues smaller than $ (n-1)^2 / 4 $.



\section{An abstract result on exponential decay} \label{theoremeabstrait}
\setcounter{equation}{0}
The purpose of this section is to prove Theorem \ref{theoreme} below which roughly asserts that $ L^2 $ solutions of $ P u = 0 $,
for a certain class of Schr\"odinger operators $ P $ on a half-line (see (\ref{formedePexemple}) below) with operator valued coefficients, decay super exponentially at infinity. We will see in Section \ref{section3} that the proof of Theorem \ref{vraitheoreme} is a consequence of Theorem \ref{theoreme}.

Let $ {\mathcal H} $ be a separable Hilbert space  with inner product $ \langle . , . \rangle_{\mathcal H} $ and norm $ ||\cdot ||_{\mathcal H} $. Let $ Q $ be a selfadjoint operator on $ {\mathcal H} $ such that 
\begin{eqnarray}
 Q \geq 1 . \label{normalisationQ}
\end{eqnarray} 
In the sequel, we fix a positive number $ R_0 > 0 $ and denote for simplicity
$$ L^2 {\mathcal H} := L^2 ((R_0,\infty),dr) \otimes {\mathcal H} \approx L^2 ((R_0,\infty),dr ; {\mathcal H}) . $$
We refer to \cite[Section II.4]{ReedSimon} for a definition of tensor products of Hilbert spaces and simply quote here that $ L^2 {\mathcal H} $ can  be viewed as the completion of $ C_0^{\infty} \big( (R_0,\infty) ; {\mathcal H} \big) $
equipped with the inner product
$$ (v,w) = \int_{R_0}^{\infty} \langle v (r) , w (r) \rangle_{\mathcal H} dr . $$
We note that the latter is still perfectly well defined on $ L^2 {\mathcal H} \cap C^0 ((R_0,\infty); {\mathcal H}) $, ie without any "almost everywhere" issue. Besides, in this paper, we shall only consider integrals of continuous functions (possibly $ {\mathcal H} $-valued). In the final application we shall take $ {\mathcal H} = L^2 (S) $ and $ Q = 1 - \Delta_S $, with $ \Delta_S $ the Laplace Beltrami operator on $ S $ associated to the metric $ \overline{g} $ in (\ref{decroissance}). 

 We  consider an unbounded operator on $ L^2 {\mathcal H} $ of the form
\begin{eqnarray}
 P = - \partial_r^2 + V_1 (r) \partial_r + V_2 (r) + e^{-2r} Q (r) , \label{formedePexemple}
\end{eqnarray}
where  $  (Q (r))_{r > R_0} $, $ (V_1 (r))_{r > R_0} $ and $ (V_2(r))_{r > R_0} $ are families of operators on $ {\mathcal H} $ satisfying the following conditions. For each $r$, $ Q (r) $ is an unbounded  operator such that, for all half integer  $ s = - \frac{1}{2}, 0 , \frac{1}{2} , 1, \ldots $ 
$$ Q(r) : \mbox{Dom}(Q^{s+1}) \rightarrow \mbox{Dom}(Q^s)  \ \ \mbox{is bounded}, $$
each domain $ \mbox{Dom}(Q^{s}) $ being equipped with the  norm $ || Q^{s} \varphi ||_{\mathcal H} $. We assume that $ Q(r) $ is $ C^2 $ with respect to $r$ in the sense that $ Q^s Q (r) Q^{-s-1} $ is $ C^2 $ in the strong sense for all $s$, and that there exists $ \tau_0 > 0 $ such that
\begin{equation}
||  Q^{s} \partial_r^j \big( Q (r) - Q \big) Q^{-s-1} ||_{{\mathcal H} \rightarrow {\mathcal H}}  \lesssim   r^{-\tau_0-j}, \qquad j = 0,1,2, \qquad r > R_0 .
\tag{A0} 
\end{equation}
 Notice that $ Q (r) $ is a perturbation of the selfadjoint operator $ Q $,  as $ r $ goes to infinity, but we do not assume that $ Q (r) $ is selfadjoint. Without loss of generality, by possibly replacing $ R_0 $ by a larger value, we may assume that
\begin{eqnarray}
||   \big( Q (r) - Q \big) Q^{-1} ||_{{\mathcal H} \rightarrow {\mathcal H}}  \leq 1/2, \qquad r > R_0 . \label{inversibilitepratique}
 \end{eqnarray}
 which implies that $ Q (r) Q^{-1} $ is invertible as an operator on $ {\mathcal H} $.
 
 We assume that $ V_1 (r) $  and $ V_2 (r) $ satisfy, for some $ \tau_1 , \tau_2 > 0 $ and all integer $ k \geq 0 $,
\begin{equation}
 || Q^k \partial_r^j V_1(r) Q^{-k}  ||_{{\mathcal H} \rightarrow {\mathcal H}}   \lesssim  r^{-\tau_1 - j}, \qquad j = 0, 1 , \qquad r > R_0  ,  \tag{A1}
\end{equation}
and
\begin{equation}
 || Q^k V_2(r) Q^{-k}  ||_{{\mathcal H} \rightarrow {\mathcal H}}   \lesssim  r^{-\tau_2}, \qquad r > R_0 ,  \tag{A2}
\end{equation}
 meaning more precisely that, for $ i = 1,2$,  $ V_{i} (r) $ preserves $ \mbox{Dom}(Q^k) $ for all $k$, that $ Q^k V_i(r)Q^{-k} $ is bounded on $ {\mathcal H} $ and is strongly $ C^{2-i} $ with the indicated decay rate with respect to $ r $. 

In the final application, we shall assume that $ \tau_0 > 1 $ and $ \tau_1 , \tau_2 > 2 $, but this will not be necessary at all steps.

\bigskip

\begin{defi} Given $v \in L^2 {\mathcal H} $, we say that $ \partial_r^2 v $ belongs to $ L^2 {\mathcal H} $ if there is a constant $ C $ such that
$$ \big| (v, \partial_r^2 \Phi)_{L^2 {\mathcal H}} \big| \leq C || \Phi ||_{L^2 {\mathcal H}} , $$
for all $ \Phi \in C_0^{\infty} ((R_0,\infty); {\mathcal H}) $. We then call $ \partial_r^2 v \in L^2 {\mathcal H} $ the unique $f \in L^2 {\mathcal H} $ such that
$$ (v, \partial_r^2 \Phi)_{L^2 {\mathcal H}} = (f,  \Phi)_{L^2 {\mathcal H}} \qquad \Phi \in C_0^{\infty} \big( (R_0,\infty); {\mathcal H} \big) $$
and set
$$ H^2 {\mathcal H} = \{ v \in L^2 {\mathcal H} \ | \ \partial_r^2 v \in L^2{\mathcal H} \}. $$
\end{defi}
Of course, as in the scalar case (ie $ {\mathcal H} = \Ca $), any $ v \in C^2 \big( (R_0,\infty);{\mathcal H} \big) $ such that $ v $ and $ \partial_r^2 v $ belong to $ L^2 {\mathcal H} $ will satisfy this definition.

For further use and the reader's convenience, we record without proof the following properties which are the $ {\mathcal H} $-valued analogues of the standard properties of  Sobolev spaces of scalar functions on the half line $  (R_0,\infty) $.
\begin{prop} \label{algebre} \begin{enumerate} \item{ $ H^2 {\mathcal H} \subset C^1 ((R_0,\infty);{\mathcal H}) $. }
\item{If $v \in  H^2 {\mathcal H} $, then $ \partial_r v \in L^2 {\mathcal H} $ and 
\begin{eqnarray}
 ||v(r) ||_{\mathcal H} \rightarrow 0 \qquad \mbox{and} \qquad || \partial_r v(r) ||_{\mathcal H} \rightarrow 0   \qquad \mbox{as} \ \ r \rightarrow \infty . \nonumber
 \end{eqnarray}}
\item{If $ v , w $ belong to $ H^2 {\mathcal H} $ and vanish near $ R_0 $ then
$$ (\partial_r^2 v, w)_{L^2 {\mathcal H}} = - \int_{R_0}^{\infty} \langle \partial_r v (r) , \partial_r w (r) \rangle_{\mathcal H} dr . $$}
\item{$ H^2 {\mathcal H}$ is stable by multiplication by $ C^2 $ functions of $r$ which are bounded together with their derivatives.}
\item{If $ v \in H^2 {\mathcal H} $ and if there is $ f \in C^0 \big( (R_0,\infty); {\mathcal H} \big) \cap L^2 {\mathcal H} $ such that $ (v,\partial_r^2 \Phi) = (f,\Phi) $ for all $ \Phi \in C_0^{\infty} \big( (R_0,\infty); {\mathcal H} \big) $ then
$$ v \in C^2 \big( (R_0,\infty); {\mathcal H} \big) \qquad \mbox{and} \qquad \partial_r^2 v = f . $$}
\item{Assume that  $$ v \in H^2 {\mathcal H} \cap C^0 ((R_0,\infty); \emph{Dom}(Q))  \qquad \mbox{and} \qquad e^{-2r}Q v \in L^2 {\mathcal H} $$ 
and $$ f \in C^0 ((R_0,\infty);{\mathcal H}) \cap L^2 {\mathcal H} . $$ Then
$$  e^{-2r} Q (r) \in L^2 {\mathcal H} . $$
Furthermore, we have equivalence between the fact that
\begin{eqnarray}
 (v, \partial_r^2 \Phi)_{L^2 {\mathcal H}}  = (V_1 (r)\partial_r v + V_2 (r) v  + e^{-2r}Q(r)v, \Phi)_{L^2 {\mathcal H}} + (f,\Phi)_{L^2 {\mathcal H}}, \label{Pufw}
\end{eqnarray}
for all $ \Phi \in C_0^{\infty} \big( (R_0,\infty); {\mathcal H} \big) $ and the fact that $ v \in C^2 ((R_0,\infty);{\mathcal H}) $ with
\begin{eqnarray}
\partial_r^2 v = V_1 (r)\partial_r v + V_2 (r) v + e^{-2r}Q(r)v +  f \label{Pufs}
\end{eqnarray}
holding in the usual sense.}
\end{enumerate}
\end{prop}

In item 6, that $ e^{-2r} Q (r)v $ belongs to $ L^2 {\mathcal H} $ is an easy consequence of (\ref{inversibilitepratique}). The second part of item 6 means of course that if $ P v = f $ in the distributions sense with a smooth enough $f$, then it holds in the strong sense.

\begin{defi} The two equivalent properties in item 6 will be denoted by
$$ P v = f . $$
\end{defi}


\bigskip

To state our result, we introduce the following convenient notation. Given $ v \in C^0 ((R_0,\infty);{\mathcal H}) $, 
$$ v \in L^{\infty}_{\rm exp} {\mathcal H} \qquad \stackrel{\mathrm{def}}{\Longleftrightarrow} \qquad \mbox{for all} \ N \geq 0 , \qquad \sup_{r \geq R_0 +1} e^{Nr} ||v(r)||_{\mathcal H} < \infty , $$
and, for $ 1 \leq p < \infty $ real,
$$  v \in L^p_{\rm exp} {\mathcal H} \qquad \stackrel{\mathrm{def}}{\Longleftrightarrow} \qquad \mbox{for all} \ N \geq 0 , \qquad  \int_{R_0 + 1}^{\infty} e^{Nr} ||v(r)||_{\mathcal H}^p dr  < \infty . $$
In both cases, the threshold $ R_0 +1 $ is completely irrevelant. We only mean that these are conditions at infinity.
It is also convenient to introduce the following space of families of operators on $ {\mathcal H} $,
\begin{eqnarray}
 {\mathcal B} = \left\{ B =(B(r))_{r \geq R_0+1} \ | \ \mbox{for all} \ j \geq 0, \ \  || Q^j B (r) Q^{-j} ||_{{\mathcal H} \rightarrow {\mathcal H}} \leq C_j , \ \ r \geq R_0 +1\right\} , 
 \label{definiB}
\end{eqnarray}
where, as before, we mean more precisely that $ B (r) $ preserves the domain of $ Q^j $ for all $j$ and $ r \mapsto Q^{j} B (r) Q^{-j} $ is strongly continuous.
These spaces have the following straightforward properties.
\begin{prop} \label{comparaisonexpo} For all $ p \geq 1 $, \begin{enumerate} \item{$ L^{\infty}_{\rm exp} {\mathcal H} \subset L^p_{\rm exp} {\mathcal H} \subset L^1_{\rm exp} {\mathcal H} $.}
\item{If $ v \in L^{1}_{\rm exp} {\mathcal H} $, then $ r \mapsto \int_r^{\infty} v (t)dt $ belongs to $ L^{\infty}_{\rm exp} {\mathcal H} $.}
\item{If $ v \in L^p_{\rm exp} {\mathcal H} $ and $ B \in {\mathcal B} $, then $ r \mapsto B (r) v (r) $ belongs to $ L^p_{\rm exp} {\mathcal H} $.}
\end{enumerate}
\end{prop}

Our main result is the following one.

\begin{theo} \label{theoreme} Assume in (A0), (A1) and (A2) that 
\begin{eqnarray}
 \tau_0 > 1 , \qquad \tau_1 > 2, \qquad \tau_2 > 2 . \label{decaystrong0}
\end{eqnarray} 
If 
\begin{eqnarray}
u \in  \bigcap_{k \in \Na} C^1 \big( (R_0,\infty) ; \emph{Dom}(Q^k) \big) , \label{hypothesereguliere}
\end{eqnarray}
 satisfies
\begin{eqnarray}
u  \in H^2 {\mathcal H} \qquad \mbox{and} \qquad e^{-2r} Q u \in L^2 {\mathcal H}, \label{conditiondedomaine}
\end{eqnarray}
and 
\begin{eqnarray}
 P u = 0 . \label{equationvp}
\end{eqnarray}
Then
\begin{eqnarray}
 u \in L^{\infty}_{\rm exp} {\mathcal H}, \qquad \partial_r u \in L^{\infty}_{\rm exp} {\mathcal H} , \qquad Q u \in L^{2}_{\rm exp} {\mathcal H} . \label{conditiondomaineu}
\end{eqnarray} 
\end{theo}

The condition (\ref{hypothesereguliere})  is only a regularity condition which, in concrete examples, will typically follow  from elliptic regularity, but which gives no information on the (square) integrability of $ || Q^k u (r) ||_{\mathcal H} $. The only $ L^2 $ information of this form at our disposal is given by  (\ref{conditiondedomaine}). Note also that the assumption (\ref{hypothesereguliere}) does not follow from the fact that $ u $ belongs to $ H^2 {\mathcal H} $ for this only implies that $ u \in C^1 \big( (R_0,\infty) ; \mbox{Dom}(Q^k) \big) $ with $ k = 0 $. On the other hand, the equation (\ref{equationvp}), together with the continuity of
$ r \mapsto e^{-2r}Q(r)u (r) $ imply that $u$ belongs to $ C^2 \big( (R_0,\infty); {\mathcal H} \big) $. We shall see later that $u$ actually  belongs to $ C^2 \big( (R_0,\infty) ; \mbox{Dom}(Q^k) \big)  $ for all $k$.

\bigskip
In the sequel we shall use the notation
$$ ||v||_{\mathcal E}^2 : = || \partial_r v ||^2_{L^2 {\mathcal H}} + ||e^{-r}Q^{1/2} v ||^2_{L^2 {\mathcal H}} . $$

 Let us introduce $ \varrho \in C^{\infty} (\Ra;\Ra) $ such that $ \mbox{supp} (\varrho) \subset [ 1 , \infty ) $ and $ \varrho \equiv 1 $ near infinity, and  set
$$ \varrho_R (r) := \varrho (r/R) . $$
The  following  elementary  lemma is a property of the operator $ P $ at infinity.

\begin{lemm} \label{minorationlemme} Assume that
\begin{eqnarray}
 \tau_0 > 0, \qquad \tau_1 > 1, \qquad \tau_2 > 2 . \label{decayloc1}
\end{eqnarray}
 For all $ R \gg 1 $ 
and all $ v \in  H^2{\mathcal H} $ such that
\begin{eqnarray}
  v  \in C^0 \big( (R_0,\infty) ; \emph{Dom}(Q) \big)  \qquad \mbox{and} \qquad e^{-2r}Q v \in L^2 {\mathcal H} , 
  \label{conditionborneinf}
\end{eqnarray}
we have
\begin{eqnarray}
 \emph{Re} \big( P (\varrho_R  v ), \varrho_R  v \big)_{L^2 {\mathcal H}} \geq \frac{1}{2} || \varrho_R v ||_{\mathcal E}^2
 . \label{estimeelemme15}
\end{eqnarray}
\end{lemm}

\noindent {\it Proof.} Without loss of generality, (\ref{decayloc1}) allows to assume that, for some $ \epsilon > 0 $ and all $ i = 0, 1,2 $ $ \tau_i = i + \epsilon $.
 By (A2), we have
\begin{eqnarray}
 \big| \big( V_2(r) \varrho_R v , \varrho_R v \big)_{L^2 {\mathcal H}} \big| \leq C R^{-\epsilon} || r^{-1} \varrho_R v ||_{L^2 {\mathcal H}}^2 , \label{BHardy}
\end{eqnarray}
and, by (A1),
\begin{eqnarray}
 \big| \big( V_1(r) \partial_r ( \varrho_R v ) , \varrho_R v \big)_{L^2 {\mathcal H}} \big| \leq C R^{-\epsilon} || \partial_r (\varrho_R v) ||_{L^2 {\mathcal H}} || r^{-1} \varrho_R v ||_{L^2 {\mathcal H}}^2  .\label{BHardybis}
\end{eqnarray}
By (A0) with $s = - 1/2 $, we also have
\begin{eqnarray}
 \left| \langle Q (r) v (r) , v (r) \rangle_{\mathcal H}  - \langle Q v (r) , v (r) \rangle_{\mathcal H} \right| \leq C r^{-\epsilon} || Q^{1/2} v (r) ||^{2}_{\mathcal H} . \label{CHardy}
\end{eqnarray}
Then, using  the following Hardy inequality, 
\begin{eqnarray}
 || r^{-1} \varrho_R v ||_{L^2 {\mathcal H}} \leq  2 || \partial_r \big( \varrho_R v \big) ||_{L^2 {\mathcal H}} , \label{Hardy}
\end{eqnarray}
 (\ref{BHardy}) and (\ref{BHardybis}) imply that
\begin{eqnarray*} 
  \big| \big( V_1 (r) \partial_r ( \varrho_R  v ), \varrho_R v \big)_{L^2 {\mathcal H}} \big| + \big| \big( V_2(r) \varrho_R v , \varrho_R v \big)_{L^2 {\mathcal H}} \big|
 & \leq & C R^{- \epsilon} || \partial_r (\varrho_R v ) ||_{L^2 {\mathcal H}}^2 \\ & \leq & C R^{- \epsilon} || \varrho_R v ||_{\mathcal E}^2 .
\end{eqnarray*}
On the other hand, multiplying (\ref{CHardy}) by $ e^{-2r} \varrho_R (r)^2 $ and integrating in $r$ yields
\begin{eqnarray*}
  \left| \int \langle e^{-2r} Q (r) \varrho_R (r) v (r) , \varrho_R (r) v (r) \rangle_{\mathcal H} dr  - ||e^{-r} Q^{1/2} \varrho_R v ||_{L^2{\mathcal H}}^2 \right| 
  & \leq & C R^{-\epsilon} || e^{-r} Q^{1/2} \varrho_R v  ||^{2}_{L^2 {\mathcal H}}  \\
  & \leq & C R^{- \epsilon} || \varrho_R v ||^2_{\mathcal E} .
\end{eqnarray*}
Using item 3 in Proposition \ref{algebre}, we have $ (- \partial_r^2 \varrho_R v , \varrho_R v)_{L^2 {\mathcal H}} = || \partial_r (\varrho_R v) ||^2_{L^2 {\mathcal H}} $ so we obtain
\begin{eqnarray*}
 \left| \big( P (\varrho_R  v ), \varrho_R (r) v \big)_{L^2 {\mathcal H}} - || \partial_r (\varrho_R v) ||^{2}_{L^2{\mathcal H}} - || e^{-r}Q^{1/2} (\varrho_R v) ||_{L^2 {\mathcal H}}^2 \right| \leq C R^{- \epsilon} || \varrho_R v ||_{\mathcal E}^2 ,
 \end{eqnarray*}
from which the result follows. \finpreuve

\bigskip

\begin{prop} \label{poidsclef} Let $u$ satisfy (\ref{hypothesereguliere}), (\ref{conditiondedomaine}) and (\ref{equationvp}) and assume that
\begin{eqnarray}
 \tau_0 > 0, \qquad \tau_1 > 1, \qquad \tau_2 > 2 . \label{decayopt2}
\end{eqnarray}
 Then
 for  $ R \gg 1 $  and  all integer $ k \geq 0 $, we have
$$ || \varrho_R Q^k u ||_{\mathcal E} 
< \infty . $$
\end{prop}

\noindent {\it Proof.} 
By (\ref{decayopt2}), we assume as in the previous lemma that $ \tau_i = i + \epsilon $ for some $ \epsilon > 0 $.
Fix $ k \geq 0 $ and define
$$ W_h = (1 + h Q)^{-k} Q^k , \qquad h \in (0,1] , $$
which is a family of bounded operators on $ {\mathcal H} $ preserving the domains of all powers of $ Q $ and converging in the strong sense to $ Q^k $ on $ \mbox{Dom}(Q^k) $ as $h \rightarrow 0$. For simplicity we denote
\begin{eqnarray}
 u_{R,h} := \varrho_R W_h u  . \label{defurh}
\end{eqnarray}
 We shall prove that, for all $R$ large enough,
\begin{eqnarray}
\sup_{ h \in (0,1]} ||u_{R,h}||_{\mathcal E} < \infty . \label{borneapriori}
\end{eqnarray}
For future reference, we already record that $ u_{R,h} $ satisfies 
\begin{eqnarray}
u_{R,h} \in C^0 ((R_0,\infty); \mbox{Dom}(Q)), \qquad e^{-2r} Q  u_{R,h} \in L^2 {\mathcal H} \label{detailenplus}
\end{eqnarray}
since $ W_h $ is bounded and preserves $ \mbox{Dom}(Q) $, and $ e^ {-2r} Q u_{R,h} = \varrho_R (r) W_h e^{-2r} Q u $ belongs to $ L^2 {\mathcal H} $ by (\ref{conditiondedomaine}).
 We observe next that, uniformly with respect to $h$,
\begin{eqnarray}
 || Q^{-1/2} W_h \big( Q (r)-Q \big) W_h^{-1} Q^{-1/2} ||_{{\mathcal H} \rightarrow {\mathcal H}} \leq C r^{-\epsilon}  . \label{premierebornelisse}
\end{eqnarray}
The latter follows from (A0), with $s$ half integer, and the easily verified fact that 
$$ W_h \big( Q (r)-Q \big) W_h^{-1} = \sum_{m=0}^k C_k^m \frac{h^m Q^m}{(1+hQ)^k } Q^{k-m} \left( Q (r) - Q \right)  Q^{m-k} ,  $$
where $ (hQ)^m / (1+hQ)^k $ is bounded, uniformly in $h$.
Similarly, using (A2) we also have
\begin{eqnarray}
 || W_h V_2(r) W_h^{-1} ||_{{\mathcal H} \rightarrow {\mathcal H}} \leq C r^{-2 - \epsilon}  , \label{deuxiemebornelisse}
\end{eqnarray}
and using (A1)
\begin{eqnarray}
 || W_h V_1(r) W_h^{-1} ||_{{\mathcal H} \rightarrow {\mathcal H}} \leq C r^{-1 - \epsilon}  . \label{troisiemebornelisse}
\end{eqnarray}
Since $ Q $ and $ W_h $ commute, we see that
\begin{eqnarray} 
 W_h Q (r) W_h^{-1} & = & Q (r) + (Q-Q(r)) + W_h (Q(r)-Q) W_h^{-1} \nonumber  \\
 & = : & Q (r) + \Psi_h (r) ,  \label{proximiteini1}
\end{eqnarray}
where, by (A0) and (\ref{premierebornelisse}),
\begin{eqnarray}
|| Q^{-1/2} \Psi_h (r) Q^{-1/2} ||_{{\mathcal H} \rightarrow {\mathcal H}} \leq C r^{-\epsilon}  . \label{bornepourastuce1}
\end{eqnarray}
Since $ \partial_r $ and $ W_h $ commute, we also see that
\begin{eqnarray}
 W_h \left( - \partial_r^2 + V_1 (r) \partial_r + V_2 (r) \right) W_h^{-1} = - \partial_r^2 + V_1 (r) \partial_r + V_2 (r) + V_{1,h}(r) \partial_r + V_{2,h}(r) ,
 \label{proximiteini2}
\end{eqnarray} 
with
$$ V_{i,h} (r)  := W_h V_i (r) W_h^{-1} - V_i (r) , \qquad i = 1,2 , $$
satisfying uniformly with respect to $h$ 
\begin{eqnarray}
 ||  V_{i,h} (r) ||_{{\mathcal H} \rightarrow {\mathcal H}} \leq C r^{-i-\epsilon} , \label{bornepourastuce2}
\end{eqnarray}
by  (A1), (A2), (\ref{deuxiemebornelisse}) and (\ref{troisiemebornelisse}).
Therefore, using (\ref{proximiteini1}) and (\ref{proximiteini2}), we can write
$$ P =  W_h P W_h^{-1}  + e^{-2r} \Psi_h (r) + V_{1,h} (r) \partial_r + V_{2,h} (r) , $$
so using (\ref{bornepourastuce1}) and the Hardy inequality (\ref{Hardy}) together with (\ref{bornepourastuce2}), we have
\begin{eqnarray}
\left| \left( P  u_{R,h} ,  u_{R,h}  \right)_{L^2 {\mathcal H}} - \left( W_h P W_h^{-1}  u_{R,h} ,  u_{R,h}  \right)_{L^2 {\mathcal H}} \right| & \leq  & C R^{-\epsilon} ||u_{R,h} ||_{\mathcal E}^2 , \nonumber \\
 & \leq  & \frac{1}{4} ||u_{R,h} ||_{\mathcal E}^2 , \label{righthandsidequart}
\end{eqnarray}
by choosing $R$ large enough independently of $h$.
%
%
%
%
%
 On the other hand, the equation (\ref{equationvp}) and (\ref{defurh}) imply that
\begin{eqnarray}
 W_h P W_h^{-1} u_{R,h} & = & W_h P \big( \varrho_R u \big) \nonumber \\
 & = & W_h \left( - 2 \varrho^{\prime}_R \partial_r  u - \varrho_R^{\prime \prime}  u + \varrho^{\prime}_R  V_1 (r) u \right)  \nonumber \\
 & = &  - 2 \varrho^{\prime}_R  \partial_r W_h   u - \varrho_R^{\prime \prime} W_h u +   \big( W_h V_1 (r) W_h^{-1} \big) \varrho^{\prime}_R W_h u . \label{avecchamp}
\end{eqnarray}
 Using the uniform bound (\ref{troisiemebornelisse}), (\ref{avecchamp})  implies that
we can find a constant $C $ independent of $h$ such that
\begin{eqnarray}
 \left| \left( W_h P W_h^{-1}  u_{R,h} ,  u_{R,h}  \right)_{L^2 {\mathcal H}} \right| \leq C \int_{{\rm supp} \ \varrho_R^{\prime}}
  || \partial_r Q^k u (r) ||_{\mathcal H}^2 + ||Q^k u (r) ||^2_{\mathcal H} dr \qquad h \in (0,1] . \label{borneunif2}
\end{eqnarray}
Using (\ref{hypothesereguliere}) and the compact support of $ \varrho_R^{\prime} $, the right hand side of  (\ref{borneunif2}) is finite, hence 
(\ref{righthandsidequart}) implies that for some constant $C_R$ independent of $h$,
\begin{eqnarray}
\big| \left( P u_{R,h} , u_{R,h} \right)_{L^2 {\mathcal H}} \big| \leq \frac{1}{4} || u_{R,h} ||_{\mathcal E}^2 + C_R, \qquad h \in (0,1] . \label{borneuniformeh}
\end{eqnarray}
If we assume additionally that $R$ is large enough so that (\ref{estimeelemme15}) holds true, then (\ref{estimeelemme15}) applied to $ u_{R,h} $ (recall that $ u_{R,h} $ satisfies the conditions (\ref{conditionborneinf}) by (\ref{detailenplus})) and (\ref{borneuniformeh}) imply that $ || u_{R,h} ||_{\mathcal E}^2 \leq 4 C_R $ hence (\ref{borneapriori}).
The result  is then a straightforward consequence of the Fatou Lemma. \finpreuve

\bigskip

A first very useful consequence of Proposition \ref{poidsclef} is the following.
\begin{prop} \label{localisationfrequence} Let $ \tau_0,\tau_1,\tau_2 $ and $u$ be as in Proposition \ref{poidsclef}.
\begin{itemize}
\item{ For all $ k \geq 0 $ and all $ r > R_0 $,
\begin{eqnarray}
 \int_r^{\infty} || e^{-2t}Q^k u(t) ||_{\mathcal H} dt < \infty , \label{pratiqueL11} \\
 \int_r^{\infty} || e^{-2t}Q^k \partial_t u(t) ||_{\mathcal H} dt < \infty , \label{pratiqueL12}
\end{eqnarray}}
\item{ For all $ f \in C^0(\Ra) \cap L^{\infty}(\Ra) $, all $ B \in {\mathcal B} $, and for all integers $ k \geq 0 $ and $ l = 0,1 $, we have
\begin{eqnarray}
  f = 0 \ \ \mbox{near} \ \ 0 \qquad \Longrightarrow \qquad  f (e^{-2r}Q) B (r) Q^k \partial^l_r u \in L^2_{\rm exp} {\mathcal H} . \label{microlocalisation}
\end{eqnarray} }
\end{itemize}
\end{prop}

The property (\ref{microlocalisation}) means  that, up to exponentially decaying terms, the solution to $ P u = 0 $ is spectrally localized where $ e^{-2r} Q $ is small. 

\bigskip

\noindent {\it Proof.} Let us consider (\ref{pratiqueL11}). By (\ref{hypothesereguliere}), it suffices to show that (\ref{pratiqueL11}) holds for all $r$ large enough or, equivalently, to check that the estimate holds with $ \varrho_R u $ instead of $u$. By (\ref{normalisationQ}) we have $ ||Q^{-1/2}||_{{\mathcal H} \rightarrow {\mathcal H}} \leq 1 $ hence
$$ ||e^{-2t} \varrho_R(t) Q^k  u (t)||_{\mathcal H} \leq e^{-t} || e^{-t} \varrho_R (t)Q^{k+1/2} u (t) ||_{\mathcal H} , $$
so the conclusion follows from Proposition \ref{poidsclef} and the Cauchy-Schwarz inequality (with respect to $t$). The proof of (\ref{pratiqueL12}) is analogous.

Let us now prove (\ref{microlocalisation}). This is again a condition at infinity, so we may replace $ u $ by $ \varrho_R u $ in the estimate. Note also that $ r \mapsto f (e^{-2r} Q) $ is strongly continuous, so  $ r \mapsto f (e^{-2r}Q)B(r) Q^k \partial_r^l u (r) $ is continuous by (\ref{hypothesereguliere}).
Assume first that $ l = 0 $. Write, for any arbitrary $ N \geq 0 $, 
$$ f (e^{-2rQ}) = f_N (e^{-2r}Q) e^{-2Nr}Q^N, \qquad \mbox{with} \qquad f_N (r,\mu) = \frac{f (e^{-2r}\mu)}{  (e^{-2r} \mu )^N  } , $$
and  observe that $ ||f_N(r,\mu)|| \leq C_N $ for all $r > R_0$ and all $ \mu \geq 1 $. Then (\ref{microlocalisation}) follows by writing 
\begin{eqnarray}
 e^{(2N-1)r} f(e^{-2r}Q) Q^k \varrho_R (r) u (r) =  f_N (r,Q) B_N (r) \left( e^{-r} Q^{N+k+\frac{1}{2}} \varrho_R(r) u (r) \right) , \label{droitesepare}
\end{eqnarray}
with
$$  B_N (r) = Q^N B (r) Q^{-N- \frac{1}{2}} $$
and by using that the bracket in right hand side of (\ref{droitesepare}) belongs to $ L^2 {\mathcal H} $ by Proposition \ref{poidsclef}. If $ l = 1 $, we proceed similarly by writing
$$ e^{2Nr} f(e^{-2r}Q) B(r) Q^k \partial_r \left(\varrho_R (r) u (r) \right) = f_N (r,Q) \widetilde{B}_N (r) Q^{k+N} \partial_r \left( \varrho_R (r) u (r) \right) , $$
with $ \widetilde{B}_N (r) = Q^N B(r)Q^{-N} $. This completes the proof.
 \finpreuve

\bigskip

The following proposition will be very useful in the sequel. To state it, we introduce the notation
$$ v \stackrel{L^2_{\rm exp} {\mathcal H}}{\equiv} w \qquad \Longleftrightarrow \qquad v - w \in L^2_{\rm exp} {\mathcal H} . $$
We also recall that $ {\mathcal B} $ is defined in (\ref{definiB}).
\begin{prop} \label{offbottom} Let $ \tau_0,\tau_1,\tau_2 $ and $ u $ be as Proposition \ref{poidsclef}. Fix 
$$ B_1 , B_2 \in {\mathcal B}   
\qquad \mbox{and} \qquad \chi \in C_0^{\infty}, \qquad \chi = 1  \ \mbox{ near } \ 0 . $$ Then, for all  $ k \geq 1 $ and $ j \geq 0$, we have
\begin{eqnarray}
B_1(r) \int_r^{\infty} B_2 (t) e^{-2kt} Q^k u (t)dt \in \emph{Dom}(Q^j) , \label{verifdomaine}
\end{eqnarray}
and
\begin{eqnarray}
Q^j B_1(r)\int_r^{\infty} B_2 (t) e^{-2kt} Q^k u (t) dt & \stackrel{L^2_{\rm exp} {\mathcal H}}{\equiv} & Q^j B_1 (r)
\int_r^{\infty} \chi (e^{-2t}Q) B_2 (t) e^{-2kt} Q^k   u (t) dt ,
\label{sansj} \\
  & \stackrel{L^2_{\rm exp} {\mathcal H}}{\equiv} &   \chi (e^{-2r}Q) Q^j B_1 (r) 
 \int_r^{\infty} B_2 (t) e^{-2kt} Q^k  u (t) dt , \label{avecj} \\
  & \stackrel{L^2_{\rm exp} {\mathcal H}}{\equiv} &   \chi (e^{-2r}Q) Q^j B_1 (r) \int_r^{\infty} \chi (e^{-2t}Q) B_2 (t) e^{-2kt} Q^k  u (t) dt  . \nonumber
 \end{eqnarray}
\end{prop}

This proposition means that, up to super exponentially decaying terms, integrals as in (\ref{verifdomaine}) can be spectrally localized where $ e^{-2r}Q $ (and/or $ e^{-2t}Q $) is small.

\bigskip

\noindent {\it Proof.} Let us show (\ref{verifdomaine}). It suffices to observe that, for any $ \varphi \in \mbox{Dom}(Q^j) $,
\begin{eqnarray*}
\left| \left\langle B_1(r)\int_r^{\infty} B_2(t)e^{-2kt} Q^k u (t)dt,Q^j \varphi  \right\rangle_{\mathcal H} \right| & =  & \left| \int_r^{\infty} e^{-2kt} \left\langle B_1(r)B_2(t) Q^k u (t)dt,Q^j \varphi  \right\rangle_{\mathcal H} dt \right|, \\
& = & \left| \int_r^{\infty} e^{-2kt} \left\langle B_j (r,t) Q^{k+j} u (t)dt, \varphi  \right\rangle_{\mathcal H} dt \right|, \\
& \leq & C || \varphi ||_{\mathcal H}
\end{eqnarray*}
with
$$ B_j (r,t) = Q^j B_1 (r) B_2 (t) Q^{-j} , $$
using the Cauchy-Schwarz inequality, (\ref{pratiqueL11}) and the fact that $ || B_j (r,t) ||_{{\mathcal H} \rightarrow {\mathcal H}} $ is bounded with respect to $r$ and $t$, which follows from definition of ${\mathcal B}$ (see (\ref{definiB})).
 We next prove only (\ref{sansj}) and (\ref{avecj}) since they imply easily the last identity. The difference between the two sides of (\ref{sansj}) reads
\begin{eqnarray*}
I_{k,j}(r) :=  B_{1,j} (r) \int_r^{\infty} (1- \chi (e^{-2t}Q))  B_{2,j} (t)  e^{-2kt}   Q^{k+j} u (t) dt .
\end{eqnarray*}
with
$$ B_{1,j}(r) =  Q^j B_1 (r) Q^{-j}, \qquad B_{2,j} (t) = Q^j B_2 (t) Q^{-j} , $$
which both belong to $ {\mathcal B} $.
By  the Cauchy Schwartz inequality, we obtain for all $ N > 0$,
$$ || I_{k,j}(r)||_{\mathcal H} \leq C_{j,N} \left( \int_r^{\infty} e^{-2Nt-2kt} dt  \right)^{1/2} \left( \int_r^{\infty} e^{2Nt}||(1- \chi (e^{-2t}Q)) B_{2,j}(t) Q^{k+j} u (t)||^2_{\mathcal H} dt \right)^{1/2} . $$
By (\ref{microlocalisation})  the second integral in the right hand side is bounded with respect to $r$, so we obtain the bound  $ || I_{j,k}(r) ||_{\mathcal H} \lesssim e^{-Nr} $, which shows that $ I_{k,j} $ belongs to $ L^{\infty}_{\rm exp} {\mathcal H} $ hence to $ L^2_{\rm exp} {\mathcal H} $.
By setting, 
$$ f (\lambda) = \frac{1- \chi (\lambda)}{\lambda^{N}} , $$
which is bounded, the
 difference between the two sides of (\ref{avecj}) can be written for all $ N \geq 0$,
\begin{eqnarray}
 II_{j,k}(r) & := & (e^{-2 r} Q )^N  f (e^{-2r}Q) Q^j B_{1}(r) \int_r^{\infty} B_{2} (t)  e^{-2kt} Q^{k} u (t) dt , \nonumber \\
 & = & e^{-2 N r} f (e^{-2r}Q) B_{1,j+N}(r) \int_r^{\infty} B_{2,j+N} (t)  e^{-2kt} Q^{N+j+k} u (t) dt , \label{IIjk}
\end{eqnarray} 
 Proceeding as in the case of $ I_{k,j} $, we see that  the $ {\mathcal H} $ norm of the integral in (\ref{IIjk}) is bounded with respect to $r$ which proves that
$ II_{j,k} $ belongs to $ L^{\infty}_{\rm exp} {\mathcal H} $ hence to $ L^2_{\rm exp}{\mathcal H} $ .
 \finpreuve

\bigskip

To justify integrations by part and get rid of boundary terms at infinity, we shall need the following lemma.
\begin{lemm} \label{termesdebordderives}  Let $ \tau_0,\tau_1,\tau_2 $ and $u$ be as in Proposition \ref{poidsclef}. Then, for all integers $ k \geq 1 $ and $ l = 0,1 $, we have
$$ \big| \big| e^{-2r} Q^k \partial_r^l u(r) \big| \big|_{\mathcal H} \rightarrow 0 , \qquad r \rightarrow + \infty . $$
\end{lemm}

\noindent {\it Proof.} It is  sufficient to show that  the $ {\mathcal H} $-valued function of $r$
\begin{eqnarray}
 \frac{\partial}{\partial r} \left( e^{-2r}Q^k  \partial^l u(r) \right)  , \label{doitetreL1}
\end{eqnarray}
is integrable  at infinity (in the Riemann sense), for this will imply that $  e^{-2r}Q^k \partial_r^l u(r)  $ has a limit as $r$ goes to infinity which is necessarily zero by (\ref{pratiqueL11})-(\ref{pratiqueL12}). Notice that (\ref{doitetreL1}) makes sense for $l=0$ by (\ref{hypothesereguliere}) and, for $ l =1 $, by item 6 of Proposition \ref{algebre}. Let us prove the integrability of (\ref{doitetreL1}). For $l=0$, we have
$$ \frac{\partial}{\partial r} \left( e^{-2r}Q^k   u(r) \right) = -2 e^{-2r}Q^k u(r) + e^{-2r}Q^k \partial_r u (r), $$
which is integrable by (\ref{pratiqueL11}) and (\ref{pratiqueL12}). For $l=1$, the equation (\ref{equationvp}) shows that
\begin{eqnarray}
 \frac{\partial}{\partial r} \left( e^{-2r}Q^k  \partial_r u(r) \right) & = &  e^{-2r}Q^k \left( -2 \partial_r u(r) +  V_1(r) \partial_r u (r) +   V_2 (r) u(r) + e^{-2r}  Q (r) u (r) \right) \nonumber
 \\
 & = & B (r) e^{-2r}Q^k \partial_r u(r) + \widetilde{B}(r) e^{-2r}Q^k u (r) + \widehat{B}(r) e^{-4r} Q^{k+1} u (r) , \label{chech}
\end{eqnarray}
where
$$ B (r) = - 2 + Q^k V_1 (r) Q^{-k}, \qquad \widetilde{B}(r) = Q^k V_2 (r) Q^{-k}, \qquad \widehat{B}(r) = Q^k Q(r) Q^{-k-1}, $$
all belong to $ {\mathcal B} $, by (A0), (A1) and (A2). Thus (\ref{pratiqueL11}) and (\ref{pratiqueL12}) show that all terms of (\ref{chech}) are integrable which completes the proof.
 \finpreuve

\bigskip

From now on we assume that $ \tau_0, \tau_1, \tau_2 $ are chosen as in Theorem \ref{theoreme} and we fix $ \epsilon > 0 $ such that 
\begin{eqnarray}
 \tau_0 \geq 1 + \epsilon, \qquad \tau_1 \geq 2 + \epsilon , \qquad \tau_2 \geq 2 + \epsilon . \label{choixepsilon}
\end{eqnarray} 
Let us then introduce the following subspace of $ {\mathcal B} $ (defined by (\ref{definiB})) 
\begin{eqnarray*}
{\mathcal S}  =  \left\{ S =(S(r))_{r \geq R_0 + 1} \in {\mathcal B} \ | \ || Q^j \partial_r^l  S(r)  Q^{-j} ||_{{\mathcal H} \rightarrow {\mathcal H}} \leq C_j r^{-2-\epsilon}, \ \ j \geq 0, \ \ 1 \leq l \leq 2 \right\},
\end{eqnarray*} 
where we implicitly assume that, for all $j$,  $ r \mapsto Q^j S (r) Q^{-j} $ is strongly $ C^2 $. Notice that, since $ || Q^k \partial_r  S(r)  Q^{-k} ||_{{\mathcal H} \rightarrow {\mathcal H}} $ is integrable at infinity, $ Q^k S(r)Q^{-k} $ has a limit as $ r\rightarrow + \infty $. 
A useful example of element of $ {\mathcal S} $ is given by
\begin{eqnarray}
 \overline{S} (r) := Q(r) Q^{-1} , \label{fauxunitaire}
\end{eqnarray}
as can be seen easily using (A0). For $q=0,1$, it is also convenient to introduce the sets
$$ {\mathcal A}_q = \left\{ v \in C^0 \big( (R_0,\infty) ; {\mathcal H} \big) \ | \ ||v(r) ||_{\mathcal H} \leq C r^{-q-\epsilon} \sup_{t \geq r} ||u(t)||_{\mathcal H} \right\} , $$
where $u$ is the function considered in Theorem \ref{theoreme}. Notice that $ \sup_{t \geq r} ||u(t)||_{\mathcal H} $ is well defined since $ u \in H^2 {\mathcal H} $ implies that  $ t \mapsto ||u(t)||_{\mathcal H} $ is bounded, by Proposition \ref{algebre}.

The following technical proposition will be useful to prove Theorem \ref{theoreme} using an induction argument.

\begin{prop} \label{intermediaireclef} Let  $ \tau_0,\tau_1,\tau_2 $ and $u$ be as in Theorem \ref{theoreme}. Let 
$$ k \geq 1, \qquad S \in {\mathcal S}   
\qquad \mbox{and} \qquad \chi \in C_0^{\infty}, \qquad \chi = 1  \ \mbox{ near } \ 0 . $$
Define
$$ \widetilde{S}(r) = S (r) Q^k Q (r) Q^{-k-1}
. $$
 Then
\begin{eqnarray}
\big( \widetilde{S}(r) \big)_{r \geq R_0 + 1} \in {\mathcal S}, \label{statementeasy} 
\end{eqnarray}
and 
\begin{eqnarray}
(2k)^2 \int_r^{\infty} S (t) e^{-2kt} Q^k u (t) dt & = & \int_r^{\infty} \widetilde{S} (t) e^{-2(k+1)t} Q^{k+1} u (t) dt   
 - \partial_r \left( S (r) e^{-2kr} Q^k \chi (e^{-2r}Q) u (r) \right)  \nonumber \\
& & +  2 S (r) e^{-2kr} Q^k  \partial_r u (r) \ \ \ \emph{mod} \ \ \ L^2_{\rm exp} {\mathcal H} + {\mathcal A}_1 . \label{clefrecurrence}
\end{eqnarray}
\end{prop}

\noindent {\it Proof.} The statement (\ref{statementeasy}) follows easily from (A0). Let us prove (\ref{clefrecurrence}). Integrating by part twice in the left hand side of (\ref{clefrecurrence}), using $ (2k)^{-2} \partial_t^2 (e^{-2kt}) = e^{-2kt} $ and 
Lemma \ref{termesdebordderives} to handle the boundary terms at infinity, we obtain
\begin{eqnarray}
(2k)^2 \int_r^{\infty} S (t) e^{-2kt} Q^k  u (t) dt & = &  \int_r^{\infty} e^{-2kt}  \left( S^{\prime \prime}(t) Q^k u (t) + 2 S^{\prime}(t)Q^k \partial_t u (t) + S (t)Q^k \partial_t^2 u(t) \right) dt \nonumber \\
& & 
 + \left( S^{\prime}(r) + 2k S (r) \right) e^{-2 k r}Q^k u (r) + S (r) e^{-2kr}Q^k \partial_r u (r) . \label{2substitutions}
\end{eqnarray}
Inside the integral in the right hand side of (\ref{2substitutions}), we replace the term $ \partial_t^2 u (t) $  by its expression given by the equation (\ref{equationvp}). This yields
\begin{eqnarray}
\int_r^{\infty} S (t) e^{-2kt} Q^k \partial_t^2 u (t) dt & = &    
  \int_r^{\infty }B_1 (t) e^{-2kt}Q^k \partial_t u (t) dt + \int_r^{\infty} B_2 (t) e^{-2kt}Q^k u (t)dt \nonumber  \\
&  & +  \int_r^{\infty} e^{-2(k+1)t}  \widetilde{S}(t) Q^{k+1} u (t)  dt , \label{1substitution}
\end{eqnarray}
with
$$ B_i (t) = S (t)Q^k V_i (t) Q^{-k}, \qquad i = 1,2 . $$
We next get rid of the terms $ \partial_t u (t) $ inside the integrals in (\ref{2substitutions}) and  (\ref{1substitution}) by integrating by part. These integrals read
\begin{eqnarray}
\int_r^{\infty} e^{-2kt}   S^{\prime}(t) Q^k \partial_t u (t) dt  =  - S^{\prime}(r) e^{-2kr}Q^k u (r) + \int_r^{\infty} \left(2k S^{\prime}(t)
-S^{\prime \prime}(t) \right) e^{-2kt} Q^k u (t)dt , \label{plusnaturel1}
\end{eqnarray}
and
\begin{eqnarray}
\int_r^{\infty }B_1 (t) e^{-2kt}Q^k \partial_t u (t) dt  =  - B_1(r) e^{-2kr}Q^k u (r) + \int_r^{\infty} \left(2k B_1(t)
-B_1^{ \prime}(t) \right) e^{-2kt} Q^k u (t)dt . \label{plusnaturel2}
\end{eqnarray}
By (\ref{choixepsilon}), we have 
$$ ||S^{\prime \prime}(r) ||_{{\mathcal H} \rightarrow {\mathcal H}} + ||S^{\prime}(r) ||_{{\mathcal H} \rightarrow {\mathcal H}} + 
||B_1^{\prime}(r) ||_{{\mathcal H} \rightarrow {\mathcal H}} + ||B_1(r) ||_{{\mathcal H} \rightarrow {\mathcal H}} + ||B_2(r) ||_{{\mathcal H} \rightarrow {\mathcal H}} \lesssim r^{-2-\epsilon} , $$
hence, using (\ref{microlocalisation}),
we see that (\ref{plusnaturel1}), (\ref{plusnaturel2}) as well as the integrals involving $ S^{\prime \prime} $ in (\ref{2substitutions}) and $ B_2 $ in (\ref{1substitution}) all belong to $ {\mathcal A}_1 + L^2_{\rm exp} {\mathcal H} $. The same holds for $ S^{\prime} (r) e^{-2kr}Q^k u (r) $ in (\ref{2substitutions}) hence (\ref{2substitutions}) reads
\begin{eqnarray}
(2k)^2 \int_r^{\infty} S (t) e^{-2kt} Q^k  u (t) dt & = &   \int_r^{\infty}   \widetilde{S}(t)e^{-2(k+1)t}Q^{k+1} u (t) dt  + 2k S (r)  e^{-2 k r}Q^k u (r) \nonumber
\\ & & \ \ + S (r) e^{-2kr}Q^k \partial_r u (r)  \qquad \mbox{mod} \qquad {\mathcal A}_1 +  L^2_{\rm exp} {\mathcal H} . \nonumber
\end{eqnarray}
The conclusion follows once observed that, 
\begin{eqnarray*}
 S (r) e^{-2kr}Q^k  u (r) & = & S (r) e^{-2kr}Q^k \chi (e^{-2r}Q)  u (r) \qquad \mbox{mod} \qquad L^2_{\rm exp} {\mathcal H} , \\
\end{eqnarray*} 
 by (\ref{microlocalisation}), and
\begin{eqnarray*}
 2k S (r)  e^{-2 k r}Q^k \chi (e^{-2r}Q) u (r) 
 & = & - \partial_r \left( S (r)  e^{-2 k r}Q^k \chi (e^{-2r}Q) u (r) \right) + S (r) e^{-2k r} Q^k \chi (e^{-2r}Q) \partial_r u (r)  \\
& & + e^{-2kr} \left( S^{\prime}(r)  \chi (e^{-2r}Q)  -2  S (r) e^{-2r}Q \chi^{\prime} (e^{-2r}Q) \right) Q^k u (r) , 
\end{eqnarray*} 
whose second line belongs to $ {\mathcal A}_1 + L^2_{\rm exp} {\mathcal H} $,
by (\ref{microlocalisation}) and the fact that $ || S^{\prime} (r) ||_{{\mathcal H} \rightarrow {\mathcal H}} \lesssim r^{-1-\epsilon} $, and where in the first line
$ e^{-2k r} Q^k \chi (e^{-2r}Q) \partial_r u (r) $ can be replaced by $ e^{-2k r} Q^k \partial_r u (r) $ thanks to  (\ref{microlocalisation}). This completes the proof. \finpreuve

\bigskip





 Theorem \ref{theoreme} will be proven thanks to the following proposition.


\begin{prop} \label{recurrence} Let  $ \tau_0,\tau_1,\tau_2 $ and $u$ be as in Theorem \ref{theoreme}.
 For all $ N \geq 1 $, there exist 
$$  B_1 , \ldots , B_{N-1} \in {\mathcal S} ,  $$
an integer $J_N$ and
$$  D_{j,k} \in {\mathcal S}, \qquad S_{j,k} \in {\mathcal S} , \qquad 1 \leq k \leq N , \ \ 1 \leq j \leq J_N ,$$
such that, for all $ \chi \in C_0^{\infty}(\Ra) $ satisfying
\begin{eqnarray}
 \chi = 1 \qquad \mbox{near} \ 0 ,  \label{satisfying}
\end{eqnarray}
we have
\begin{eqnarray}
\partial_r u (r)  &  = & \sum_{k=1}^{N} \sum_{j=1}^{J_N} D_{k,j} (r) (e^{-2r} Q)^{N-k} \int_r^{\infty} S_{k,j} (t)   (e^{-2t}Q)^k u(t)dt \ + \nonumber \\
 &   & \ \ \ \sum_{j=1}^{N-1} \partial_r \left( B_j (r) (e^{-2r} Q)^j \chi (e^{-2r}Q) u (r) \right) \ \ \ \emph{mod} \ \ \  {\mathcal A}_1 + L^2_{\rm exp} {\mathcal H}  ,
 \label{recurrenceN}
\end{eqnarray}
for all $r$ large enough.
\end{prop}


Before starting the proof of Proposition \ref{recurrence}, we record the following useful simple computation in which we use the notation (\ref{fauxunitaire}).
By integrating (\ref{equationvp}) on $[r,+\infty)$ and using (\ref{choixepsilon}) (plus item 2 of Proposition \ref{algebre} to handle the boundary terms at infinity in the integration by part),
we have
\begin{eqnarray}
 - \partial_r u (r) &= & \int_r^{\infty} e^{-2t} Q (t) u (t) dt + \int_r^{\infty} V_1 (t) \partial_t u (t) d t  + \int_r^{\infty} V_2(t)u(t)dt   \label{conclusioncomplete} \\
 & = &  \int_r^{\infty} e^{-2t} Q (t) u (t) dt + \int_r^{\infty} \big( V_2(t) - \partial_t V_1 (t) \big)u(t)dt  - V_1 (r) u (r)
\label{conclusioncompleteintermediaire}  \\
 & = & \int_r^{\infty} e^{-2t} \overline{S} (t) Q u (t) dt \qquad \mbox{mod} \qquad {\mathcal A}_1 . \label{initialisation}
\end{eqnarray}

\bigskip

\noindent {\bf Proof of Proposition \ref{recurrence}.} We proceed by induction. The case $ N = 1 $ is a direct consequence of (\ref{initialisation}).
To go from step $ N $ to $N+1 $, it is sufficient to consider  the first sum in the right hand side of (\ref{recurrenceN}).
We focus on a single term,
$$ T (r) := D (r) (e^{-2r}Q)^{N-k} \int_r^{\infty} S (t) (e^{-2t}Q)^k u (t)dt , $$
 where we drop the indices $k,j$ to simplify the notation. It suffices to show that, for some $ B_N  \in {\mathcal S} $ and $ D_1  , D_{k+1}, S_{k+1} \in {\mathcal S} $ independent of $ \chi $, we can write 
\begin{eqnarray}
T (r)  &  = &  \partial_r \left( B_N (r) (e^{-2r} Q)^N \chi (e^{-2r}Q) u (r) \right) + D_1 (r) (e^{-2r} Q)^{N} \int_r^{\infty} \overline{S} (t)   e^{-2t}Q u(t)dt \ + \nonumber \\
 &   &  D_{k+1} (r) (e^{-2r} Q)^{N-k} \int_r^{\infty} S_{k+1} (t)   \big( e^{-2t}Q \big)^{k+1} u(t)dt  \ \ \ \mbox{mod} \ \ \  {\mathcal A}_1 + L^2_{\rm exp} {\mathcal H}  .
 \label{recurrenceNbis}
\end{eqnarray}
Let us prove this fact. By Proposition \ref{offbottom}, we have
\begin{eqnarray}
 T(r)  =   D (r) (e^{-2r} Q)^{N-k} \chi (e^{-2r}Q) \int_r^{\infty} S (t)   (e^{-2t}Q)^k u(t)dt \qquad \mbox{mod} \ \ \ L^2_{\rm exp} {\mathcal H}  . \label{righthandsideTmod}
\end{eqnarray}
Using Proposition \ref{intermediaireclef} and the fact that the space $ L^2_{\rm exp} {\mathcal H} + {\mathcal A}_1 $ is preserved by the action of $ D(r) (e^{-2r}Q)^{N-k} \chi (e^{-2r}Q) $, the right hand side of (\ref{righthandsideTmod}) is a linear combination  of
\begin{eqnarray*}
T_1 & : =& D (r) (e^{-2r} Q)^{N-k} \chi (e^{-2r}Q) S (r) e^{-2kr}Q^k \partial_r u (r) ,  \\
T_2 & := &  D (r) (e^{-2r} Q)^{N-k} \chi (e^{-2r}Q) \partial_r \left( S (r) e^{-2kr} Q^k \chi (e^{-2r}Q) u (r) \right) , \\
T_3 & := &   D (r) (e^{-2r} Q)^{N-k} \chi (e^{-2r}Q) \int_r^{\infty} \widetilde{S} (t) e^{-2(k+1)t} Q^{k+1} u (t) dt ,
\end{eqnarray*}
and of a term in $ L^2_{\rm exp} {\mathcal H} + {\mathcal A}_1 $.  Here $ \widetilde{S} \in {\mathcal S} $ is defined as in Proposition \ref{intermediaireclef}.
By setting $$ \widehat{S}(r) = Q^{N-k} S (r) Q^{k-N} , $$
we have
\begin{eqnarray}
T_1 & = & D (r)  \chi (e^{-2r}Q) \widehat{S} (r) \big( e^{-2r}Q \big)^N \partial_r u (r) \nonumber \\
 & = & D (r)  \widehat{S} (r) \big( e^{-2r}Q \big)^N \partial_r u (r) \qquad \mbox{mod} \qquad L^2_{\rm exp} {\mathcal H}, \nonumber  \\
 & = &  D (r)  \widehat{S} (r) \chi (e^{-2r}Q) \big( e^{-2r}Q \big)^N \partial_r u (r) \qquad \mbox{mod} \qquad L^2_{\rm exp} {\mathcal H}, \nonumber  \\
 & = &  D (r)  \widehat{S} (r) \chi (e^{-2r}Q) \big( e^{-2r}Q \big)^N \int_r^{\infty} e^{-2t} \overline{S} (t) Q u (t) dt \qquad \mbox{mod} \qquad {\mathcal A}_1 + L^2_{\rm exp} {\mathcal H} , \nonumber \\
 & = &  D (r)  \widehat{S} (r) \big( e^{-2r}Q)^N \int_r^{\infty} e^{-2t} \overline{S} (t) Q u (t) dt \qquad \mbox{mod} \qquad {\mathcal A}_1 + L^2_{\rm exp} {\mathcal H} ,
 \label{conclusionI}
\end{eqnarray}
which is as the right hand side of (\ref{recurrenceNbis}). Here we have used (\ref{microlocalisation})  in the second and third  lines, (\ref{initialisation}) in the fourth one and Proposition \ref{offbottom} in the fifth one jointly with the fact that $ D (r)  \widehat{S} (r) \chi (e^{-2r}Q) (e^{-2r}Q)^N $ belongs to $ {\mathcal B} $ hence preserves $ {\mathcal A}_1 $. 
As a direct consequence of Proposition \ref{offbottom}, we also have 
\begin{eqnarray}
T_3 & =  &   D (r) (e^{-2r}Q)^{N-k} \int_r^{\infty} \widetilde{S} (t) e^{-2(k+1)t} Q^{k+1} u (t) dt  \qquad \mbox{mod} \qquad L^2_{\rm exp} {\mathcal H} , \label{conclusionIII}
\end{eqnarray}
which, as $T_1$, is also of the form of the right hand side of (\ref{recurrenceNbis}). It remains to prove that $ T_2 $ is of this form too. By (\ref{microlocalisation}) and the fact that $ ||S^{\prime}(r)||_{{\mathcal H} \rightarrow {\mathcal H}} \lesssim r^{-1-\epsilon} $, we have
\begin{eqnarray}
T_2  & = &   D (r) (e^{-2r} Q)^{N-k} \chi (e^{-2r}Q)  S (r) e^{-2kr} Q^k \chi (e^{-2r}Q)(\partial_r u (r) - 2 k u (r) ) \ \ \ \mbox{mod} \ \ \ {\mathcal A}_1 + L^2_{\rm exp} {\mathcal H} , \nonumber \\
& = &   D (r) \widehat{S}(r) \big( e^{-2r} Q\big)^N  (\partial_r u (r) - 2 k u (r) ) \ \ \ \mbox{mod} \ \ \ {\mathcal A}_1 + L^2_{\rm exp} {\mathcal H} . \label{numeroderivee1}
\end{eqnarray}
Therefore, if we set
$$ B (r) = D (r) \widehat{S}(r) , $$
the contribution of $ \partial_r u (r) $ in (\ref{numeroderivee1}) is given by
\begin{eqnarray}
B(r) \big( e^{-2r} Q \big)^N \partial_r  u (r)  & = & B(r)  \big( e^{-2r} Q \big)^N \chi (e^{-2r}Q) \partial_r  u (r) \qquad \mbox{mod} \qquad L^2_{\rm exp} {\mathcal H}, \nonumber \\
& = & - B(r) \big( e^{-2r} Q \big)^N \chi (e^{-2r}Q) \int_r^{\infty} \overline{S}(t) e^{-2t}Q u (t) dt \ \ \ \mbox{mod} \ \ \ {\mathcal A}_1 + L^2_{\rm exp} {\mathcal H}, \nonumber \\
& = & - B(r) \big( e^{-2r} Q \big)^N  \int_r^{\infty} \overline{S}(t) e^{-2t}Q u (t) dt \ \ \ \mbox{mod} \ \ \ {\mathcal A}_1 + L^2_{\rm exp} {\mathcal H}, \label{contributionderivee}
\end{eqnarray}
using (\ref{microlocalisation}) in the first line, (\ref{initialisation}) in the second one and Proposition \ref{offbottom} in the third one. The contribution of $ -2k u (r) $ in (\ref{numeroderivee1}) is obtained by using
\begin{eqnarray}
2N B (r) \big( e^{-2r} Q \big)^N   u (r)  & = & 2N B (r) \big( e^{-2r} Q \big)^N \chi (e^{-2r}Q)  u (r) \qquad \mbox{mod} \qquad L^2_{\rm exp} {\mathcal H} , \nonumber \\
& = & B (r) \big( e^{-2r} Q \big)^N \chi (e^{-2r}Q) \partial_r u (r) - \partial_r \left( B (r) \big( e^{-2r} Q \big)^N \chi (e^{-2r}Q)  u (r) \right) \nonumber \\
& & \ \  \ \ \ \mbox{mod} \ \ \ \ {\mathcal A}_1 + L^2_{\rm exp} {\mathcal H} . \label{listeutile}
\end{eqnarray}
where the first line follows from (\ref{microlocalisation}) and the second one from the fact that $ B $ belongs to $ {\mathcal S} $ and again (\ref{microlocalisation}).
Using (\ref{microlocalisation}), (\ref{initialisation}), (\ref{contributionderivee}) and (\ref{listeutile}), we obtain
\begin{eqnarray*}
T_2 & = & \left( \frac{k}{N} - 1 \right) B (r) \big( e^{-2r} Q \big)^N  \int_r^{\infty} \overline{S}(t) e^{-2t}Q u (t) dt + \partial_r \left( \frac{k}{N} B (r) \big( e^{-2r} Q \big)^N \chi (e^{-2r}Q)  u (r) \right) \nonumber \\
& & \ \  \ \ \ \mbox{mod} \ \ \ \ {\mathcal A}_1 + L^2_{\rm exp} {\mathcal H}
\end{eqnarray*}
which is as the right hand side of (\ref{recurrenceNbis}). The proof is complete. \finpreuve

\bigskip

\bigskip

\noindent {\bf Proof of Theorem \ref{theoreme}.}  For any $ \chi \in C_0^{\infty} $ satisfying (\ref{satisfying}) to be fixed below  and for any fixed $ N \geq 1 $, we introduce 
$$ v_{\chi,N} \in {\mathcal A}_1, \qquad s_{\chi,N} \in L^2_{\rm exp} {\mathcal H} , $$
such that 
$$ \partial_r u \ - \left( \mbox{ sums in the right hand side of } \ (\ref{recurrenceN}) \right) \ = \ v_{\chi,N} + s_{\chi,N} . $$ By item 2 of Proposition \ref{comparaisonexpo}, we have
\begin{eqnarray}
 r \mapsto \int_r^{\infty} s_{\chi,N}(t) d t \in L^{\infty}_{\rm exp} {\mathcal H} . \label{pourreference}
\end{eqnarray}
We also clearly have
$$ \tilde{v}_{N,\chi} := \left[ r \mapsto \int_r^{\infty} v_{\chi,N}(t) d t \right] \in {\mathcal A}_0 . $$
On the other hand, let us observe that, if we set $$ \widetilde{S}_{k,j}(t) = Q^{N-k} S_{k,j}(t) Q^{k-N} , $$ 
which defines a family in $ {\mathcal S} $,
then
$$ Q^{N-k} \int_r^{\infty} e^{-2kt} S_{k,j} (t) Q^k u(t)dt = \int_r^{\infty} e^{-2kt} \widetilde{S}_{k,j} (t) Q^{N} u(t)dt $$
and we have,   by  (\ref{pratiqueL11}),
$$ \left| \left| \int_r^{\infty} e^{-2kt} \widetilde{S}_{k,j} (t) Q^N u(t)dt \right| \right|_{\mathcal H} \leq C e^{-2(k-1)r} . $$
Therefore, the first line of (\ref{recurrenceN}) is $ {\mathcal O}(e^{-2(N-1)r}) $ in ${\mathcal H} $. 
Integrating (\ref{recurrenceN}) using item 2 of Proposition \ref{algebre}, Lemma \ref{termesdebordderives} and (\ref{pourreference}), we see that, for all $ \chi  $ satisfying (\ref{satisfying}), we have 
\begin{eqnarray}
 \left| \left| \left( I_{\mathcal H} - \sum_{j=1}^N  B_j (r) (e^{-2r} Q)^j \chi (e^{-2r}Q) \right) u (r)  \right| \right|_{\mathcal H} \leq
C_{N,\chi} e^{-2 (N-1) r} + || \tilde{v}_{N,\chi}(r) ||_{\mathcal H} , \qquad r \gg 1 . \label{suchachi}
\end{eqnarray} 
By choosing $ \chi $ with support close enough to $ 0 $, the norm $ || (e^{-2r}Q)^j \chi (e^{-2r}Q) ||_{{\mathcal H} \rightarrow {\mathcal H}} $ is small uniformly in $r$  since, for each $j \geq 1$,  $ \sup_{\lambda \in \Ra} | \lambda^j \chi (\lambda) | $ is as small as wish by shrinking the support of $ \chi $ to $ \{ 0 \} $. Thus,  using that $ B_1(r) , \ldots , B_N (r) $ don't depend on $ \chi $ and are uniformly bounded, we may assume that $ \chi $ is chosen so that
$$ \left| \left| \ \sum_{j=1}^N  B_j (r) (e^{-2r} Q)^j \chi (e^{-2r}Q)   \right| \right|_{{\mathcal H} \rightarrow {\mathcal H}} \leq 1/2  , $$

Using such a $ \chi $, the left hand side of (\ref{suchachi}) is bounded below by $ ||u(r)||_{\mathcal H} /2 $ so we obtain
$$ ||u(r) ||_{\mathcal H} \leq C^{\prime}_{N,\chi} e^{-2(N-1)r} + C_{N,\chi}^{\prime} r^{-\epsilon} \sup_{t \geq r} || u(t) ||_{\mathcal H} , $$
where the left hand side can be replaced by $ \sup_{t \geq r} ||u(t) ||_{\mathcal H} $ since the right hand side is a non increasing function of $r$. Then, for $r$ large enough, the second term of the right hand side can be absorbed in the left hand side and we obtain 
$$ \sup_{r \geq t} || u(t) ||_{\mathcal H} \leq 2 C^{\prime}_{N,\chi} e^{-2(N-1)r}, \qquad r \geq R_{N,\chi} . $$
Since $N$ is arbitrary, this shows that $ u $ belongs to $ L^{\infty}_{\rm exp} {\mathcal H} $. Then, by writing
$$ Q^k u = Q^k \chi (e^{-2r}Q) u + Q^k (1 - \chi (e^{-2r}Q)) u , $$
where the first term in the right hand side belongs to $ L^{\infty}_{\rm exp} {\mathcal H} $ since $u$ does, and where the second term belongs to $ L^2_{\rm exp} {\mathcal H} $ by Proposition \ref{localisationfrequence}, we obtain that $ Q^k u $ belongs to $ L^2_{\rm exp} {\mathcal H} $. Finally, using that $ u \in L^{\infty}_{\rm exp} {\mathcal H} $ and $ Q(r)u \in L^2_{\rm exp} {\mathcal H} $ (by writing $ Q(r)u = Q(r)Q^{-1} Qu $),  (\ref{conclusioncompleteintermediaire}) implies that $ \partial_r u $ belongs to $ L^{\infty}_{\rm exp} {\mathcal H} $.
\finpreuve

\section{Proof of Theorem \ref{vraitheoreme}} \label{section3}
\setcounter{equation}{0}
We start by checking that  the  equation (\ref{equationvaleurpropre})   can be reduced to $ P u = 0 $ with $ P $ as in Section \ref{theoremeabstrait} and $u$ statisfying  the conditions (\ref{hypothesereguliere}) and (\ref{conditiondedomaine}) of Theorem \ref{theoreme}.

On $ M \setminus K \approx (R_0,+\infty) \times S $, the Laplace-Beltrami operator associated to (\ref{metriquehyperbolique}) takes the form
$$ \Delta_G = \partial_r^2  + e^{-2r} \Delta_{g(r)} +  c (r,\omega) \partial_r + (n-1) \partial_r , $$
where
$$ c (r,\omega) =  \frac{1}{2} \frac{\partial_r \mbox{det}( g(r,\omega) )}{\mbox{det}( g(r,\omega) )}, \qquad \omega \in S , $$
is intrinsincally defined, for each $r$, as a function on $S$. The Laplacian $ \Delta_G $ is symmetric with respect to the Riemannian measure which is of the form
$$ d {\rm vol}_G = e^{(n-1)r} dr d {\rm vol}_{g(r)} , $$
where, for each $r$, $  d {\rm vol}_{g(r)}  $ is the Riemannian measure on $S$ associated to the metric $g(r)$. Let us set
$$ f (r,\omega)=  \left( \frac{ {\rm det}( g(r,\omega))}{{\rm det} ( \overline{g}(\omega) )} \right)^{1/2} , $$
which is again a well defined function, so that
$$ d {\rm vol}_G = e^{(n-1)r} f dr d {\rm vol}_{\overline{g}} . $$
We thus have a unitary mapping
$$  \varphi \mapsto e^{\frac{(n-1)r}{2}} f^{1/2} \varphi  $$
between $ L^2 \left( M \setminus K , d {\rm vol}_G \right)  $ and
$$  L^ 2 ((R_0,+\infty) \times S , dr d {\rm vol}_{\overline{g}}) \approx L^2 (R_0, \infty) \otimes  L^2 (S, d {\rm vol}_{\overline{g}}) .  $$
Then it is not hard to see that $ - e^{\frac{(n-1)r}{2}} f^{1/2} \Delta_G f^{-1/2} e^{-\frac{(n-1)r}{2}} $ reads $$ - \partial_r^2 - e^{-2r} f^{1/2} \Delta_{g(r)} f^{-1/2} + \frac{ (n-1)}{2} \frac{\partial_r f}{f} + \frac{1}{2} \frac{\partial_r^2 f}{f} + \frac{1}{4} \frac{(\partial_r f)^2}{f^2} + \frac{(n-1)^2}{4} . $$
This operator is symmetric with respect to $ dr d {\rm vol}_{\overline{g}} $ and of the form
$$ - e^{\frac{(n-1)r}{2}} f^{1/2} \Delta_G f^{-1/2} e^{-\frac{(n-1)r}{2}} = P_0 + \frac{(n-1)^2}{4} $$
with $ P_0 $ as in (\ref{formedePexemple}) by taking
\begin{eqnarray*}
 {\mathcal H}  =  L^2 (S, d {\rm vol}_{\overline{g}}), \qquad Q  =  1 - \Delta_{\overline{g}} , \qquad Q (r)  =  f^{1/2} \big( 1 - \Delta_{g(r)} \big) f^{-1/2} 
\end{eqnarray*}
and 
\begin{eqnarray*}
V_1 = 0, \qquad V_2  =  \frac{ (n-1)}{2} \frac{\partial_r f}{f} + \frac{1}{2} \frac{\partial_r^2 f}{f} + \frac{1}{4} \frac{(\partial_r f)^2}{f^2} -e^{-2r} .
\end{eqnarray*}
We note also in passing that
\begin{eqnarray}
 \mbox{Dom} (Q^s) = H^{2s} (S) ,  \label{echelleSobolev}
\end{eqnarray}
where $H^{\sigma}(S)$ is the $L^2 $ Sobolev space of order $ \sigma$ on $S$.
By (\ref{decroissance}), it is not hard to check that $ f - 1 $ has the same decay properties as $ g (r) - \overline{g} $ with respect to $r$ (see (\ref{decroissance})), hence (A0),(A1) and (A2) hold with
$$ \tau_1 = \tau_2 = \tau_0 + 1 . $$
If $ V $ is an admissible perturbation, one checks similarly that
$$ P := P_0 +  e^{\frac{(n-1)r}{2}} f^{1/2} V e^{-\frac{(n-1)r}{2}} f^{-1/2}  $$ 
also satisfies the conditions (A0), (A1) and (A2). By assuming that $ \tau_0 > 1 $ (as  in Theorem \ref{vraitheoreme}), we see that (\ref{decaystrong0}) is satisfied. 
Consider now $ \psi \in L^2 (M, d {\rm vol}_G) $ such that
$$ \left( - \Delta_G  + V \right) \psi = \frac{(n-1)^2}{4} \psi , $$
and 
define 
\begin{eqnarray}
 u = \left. e^{\frac{(n-1)r}{2}} f^{1/2} \psi \right|_{M \setminus K} , \label{psiu}
\end{eqnarray} 
which satisfies $ P u = 0 $.
By standard elliptic regularity, $u$ is smooth on $ (R_0 , \infty) \times S $ hence, using (\ref{echelleSobolev}), $u$ clearly satisfies (\ref{hypothesereguliere}) 
and it remains to check (\ref{conditiondedomaine}). To prove this, we observe that the principal symbol of $ P $  is of the form
$$ p =  \rho^2 + e^{-2r} \left( \sum_{j,k=1}^{n-1} g_{jk}(r,\theta) \eta_j \eta_k + {\mathcal O}(r^{-\tau_0}|\eta|^2)  \right) ,$$
with $  \sum g_{jk}(r,\theta) \eta_j \eta_k $ the principal symbol of $ - \Delta_{g(r)} $. This expression shows that $p$ takes its values in a sector $ e^{i [-\gamma,\gamma]}[0,+\infty) $ with $ \gamma $ as small as we wish by possibly taking $ R_0 $ large enough. Choose $z \in \Ca$  at positive distance from this sector and observe that 
\begin{eqnarray}
  \left(  P - z \right) u = - z u , \label{nouvelleequation} 
\end{eqnarray}
whose right hand side belongs to $ L^2 $. Fix next $ \chi $ which is supported in $ r > R_0 $ and which is equal to $1$ near infinity. By using a parametrix for $ \left(  P^* - \bar{z} \right) $ as in \cite{Bouclet0} (see the formula (2.19) and use that $\bar{z}$ is at positive distance from the set of values of $ \bar{p} $), we can find operators $ Q $ and $R$ such that
\begin{eqnarray}
 (P^* - \bar{z}) Q = \chi + R, \label{aadjoindre} 
\end{eqnarray}
and with the property  that $ (e^{-r} \partial_{\theta})^{\alpha} \partial_r^k Q^* $ and $ (e^{-r} \partial_{\theta})^{\alpha} \partial_r^k R^* $ are bounded on $ L^2 $ if $ |\alpha| + k \leq 2 $. By taking the adjoint of (\ref{aadjoindre}) and using (\ref{nouvelleequation}), we have
$$ \chi u = - \big( z Q^*  + R^* \big) u $$
from which we get that $ D_r^2 (\chi u) $ and $ e^{-2r} \Delta_{\bar{g}}(\chi u) $ belong to $ L^2 $ hence that (\ref{conditiondedomaine}) is satisfied. \finpreuve

\bigskip

\noindent {\bf Proof of Theorem \ref{vraitheoreme}.} By (\ref{conditiondomaineu}), $ u $ and $ \partial_r u $ belong to $ L^{\infty}_{\rm exp} {\mathcal H} $
and this implies that $ e^{Cr} \psi $ and $ \partial_r (e^{Cr} \psi) $ belong to $ L^2 $ (the contribution of $f$ is easily studied using (\ref{decroissance})). Then, using that
$$   - \Delta_G  + V   = \frac{(n-1)^2}{4}   + e^{-\frac{(n-1)r}{2}}f^{-1/2} P f^{1/2} e^{\frac{(n-1)r}{2}} , $$
it suffices to show that $ P (e^{Cr} u) $ belongs to $ L^2 (R_0,\infty) \otimes {\mathcal H} $ to guarantee that $ (-\Delta_G + V) (e^{Cr} \psi) $ belongs to $ L^2 $.
This follows  on one hand from the fact that $ [P , e^{Cr} ] $ is a first order differential operator involving only derivatives with respect to $r$ and with bounded coefficients and on the other hand from the fact that $ P u = 0$. 
This completes the proof of (\ref{superexponentielfinal}) hence of Theorem \ref{vraitheoreme}. \finpreuve

\bigskip


\end{document}